\documentstyle[amscd,amssymb,verbatim,12pt]{amsart}

\theoremstyle{plain}
\newtheorem{Prop}{Proposition}[section]

\newtheorem{Cor}[Prop]{Corollary}

\newtheorem{Lem}[Prop]{Lemma}

\theoremstyle{definition}
\newtheorem{Def}[Prop]{Definition}

\theoremstyle{remark}

\def\int{\mathop{\roman{int}}}

\def\1{^{-1}}

\def\MM{{\mathcal M}}

\def\RR{{\mathbb R}}

\def\dokaz{{\bf Proof. }}
\numberwithin{equation}{section}
\begin{document}
\title[
Svarc-Milnor Lemma: a proof by definition]%
   {Svarc-Milnor Lemma: a proof by definition
}

\author{N.~Brodskiy}
\address{University of Tennessee, Knoxville, TN 37996, USA}
\email{brodskiy@@math.utk.edu}

\author{J.~Dydak}
\address{University of Tennessee, Knoxville, TN 37996, USA}
\email{dydak@@math.utk.edu}
\author{A.~Mitra}
\address{University of Tennessee, Knoxville, TN 37996, USA}
\email{ajmitra@@math.utk.edu}

\date{ March 17, 2006}
\keywords{Asymptotic dimension, coarse category, Lipschitz functions, Nagata dimension}

\subjclass{ Primary: 54F45, 54C55, Secondary: 54E35, 18B30, 54D35, 54D40, 20H15}

\begin{abstract} 

The famous \v Svarc-Milnor Lemma says that a group
$G$ acting properly and cocompactly via isometries on a length space $X$
is finitely generated and
induces a quasi-isometry equivalence $g\to g\cdot x_0$ for any $x_0\in X$. 
We redefine the concept of coarseness so that the proof of the Lemma is automatic.

\end{abstract}

\maketitle

\medskip
\medskip


Geometric group theorists traditionally restrict their attention to finitely
generated groups equipped with a word metric. A typical proof of \v Svarc-Milnor Lemma
(see \cite{Roe lectures} or \cite{Bridson-HaefligerBook}, p.140) involves such metrics. Recently, the study of large scale geometry of groups was expanded to all countable groups by usage of proper, left-invariant metrics:
in \cite{Smith} such metrics were constructed and it was shown that they all induce the same coarse structure on a group (see also \cite{DS}). The point of this note is that a
proper action of a group $G$ on a space ought to be viewed as a geometric way of creating a coarse structure on $G$. That structure is not given by a proper metric but by something very similar; a pseudo-metric where only a finite set of points may be at mutual distance
$0$. From that point of view the proof of \v Svarc-Milnor Lemma is automatic and the Lemma can be summarized as follows. There are two ways of creating coarse structures on countable groups: algebraic (via word or proper metrics) and geometric (via group actions), and both ways are equivalent.

\begin{Def}\label{LSMetricDef}

A pseudo-metric $d_X$ on a set $X$ is called a {\it large-scale metric}
(or ls-metric)
if for each $x\in X$ the set $\{y\in X \mid d_X(x,y)=0\}$ is finite.

$(X,d_X)$ is called a {\it large-scale metric space}
(or an ls-metric space)
if $d_X$ is an ls-metric.
\end{Def}

\begin{Def}\label{LSProperMetricOnGroups}
An ls-metric $d_G$ on a group $G$ is {\it proper} and {\it left-invariant}
if $d_G(g,h)=d_G(f\cdot g,f\cdot h)$ for all $f,g,h\in G$
and $\{h \mid d_G(g,h) < r\}$ is finite for all $r > 0$ and all $g\in G$.
\end{Def}
Notice $G$ must be countable if it admits a proper ls-metric.

One aspect of \v Svarc-Milnor Lemma is $G$ being finitely generated.
That corresponds to $(G,d_G)$ being metrically connected, i.e. there is $M > 0$
such that any two points in $G$ can be connected by a chain of points
separated by at most $M$.

\begin{Lem} \label{MetricConnectedness}
Suppose $d_G$ is a proper and left-invariant ls-metric on $G$.
$(G,d_G)$ is metrically connected if and only if $G$ is finitely generated.
\end{Lem}
\dokaz If $G$ is generated by a finite set $F$, put $M=\max\{d_G(1_G,f)\mid f\in F\}$.
If $(G,d_G)$ is $M$-connected, put $F=B(1_G,M+1)$.
\hfill $\blacksquare$

 \begin{Def}\label{LSUniformDef}
A function $f\colon (X,d_X)\to (Y,d_Y)$ of ls-metric spaces is called  {\it large-scale uniform}
(or ls-uniform)
if for each $r > 0$ there is $s > 0$ such that $d_X(x,y)\leq r$ implies
$d_Y(f(x),f(y))\leq s$.
\par $f$ is a {\it large-scale uniform equivalence} if there is
an ls-uniform $g:Y\to X$ such that both $g\circ f$ and $f\circ g$ are within
a finite distance from the corresponding identities.
\end{Def}

\begin{Lem} \label{LSUniformForGroups}
Suppose $(G,d_G)$ and $(H,d_H)$ are two groups equipped
with proper and left-invariant ls-metrics.
A function $f\colon (G,d_G)\to (H,d_H)$ is ls-uniform if and only
if for each finite subset $F$ of $G$ there is a finite subset $E$ of $H$
such that $x^{-1}\cdot y\in F$ implies $f(x)^{-1}\cdot f(y)\in E$
for all $x,y\in G$.
\end{Lem}
\dokaz Suppose $f$ is ls-uniform and $F$ is a finite subset of $G$.
Let $r$ be larger that all $d_X(1_G,g)$, $g\in F$.
Pick $s> 0$ such that $d_G(g,h) < r$ implies $d_H(f(g),f(h)) < s$
and put $E=\{x\in H \mid d_H(1_H,x) < s\}$.
If $x^{-1}\cdot y\in F$, then $d_G(x,y) < r$. Therefore
$s > d_H(f(x),f(y))=d_H(1_H,f(x)^{-1}\cdot f(y))$ and
$f(x)^{-1}\cdot f(y)\in E$.
Conversely, if $r > 0$ put $F=\{x\in G \mid d_G(1_G,x) < r\}$
and consider $E$ so that $x^{-1}\cdot y\in F$ implies $f(x)^{-1}\cdot f(y)\in E$.
If $s$ is bigger that all $d_H(1_H,g)$, $g\in E$, then
$d_G(x,y) < r$ implies $f(x)^{-1}\cdot f(y)\in E$ and $d_H(f(x),f(y)) < s$.
\hfill $\blacksquare$

\begin{Cor} \label{LSMetricsForGroupsAreUnique}
Given two proper and left-invariant ls-metrics $d_1$ and $d_2$ on the same group $G$,
the identity $id_G:(G,d_1)\to (G,d_2)$ is a coarse equivalence.
\end{Cor}
\dokaz The choice of $E=F$ always works for $id_G$.
\hfill $\blacksquare$

We are interested in creating proper left-invariant ls-metrics on groups $G$
using actions on metric spaces $X$ via the formula $d_G(g,h)=d_X(g\cdot x_0,h\cdot x_0)$
for some $x_0\in X$. To make $d_G$ left-invariant, a practical requirement is the action occurs
via isometries. Let's characterize the situation in which $d_G$ is a proper ls-metric.

\begin{Lem} \label{CharOflsMetrics}
Suppose $G$ acts via isometries on $X$ and $x_0\in X$. If $d_G$ is defined by 
$d_G(g,h)=d_X(g\cdot x_0,h\cdot x_0)$, then $d_G$ is a proper left-invariant ls-metric on $G$ if and only if the following conditions are satisfied:
\begin{itemize}
\item[1.]  The stabilizer $\{g\in G\mid g\cdot x_0=x_0\}$ of $x_0$ is finite.
\item[2.] $G\cdot x_0$ is topologically discrete.
\item[3.] Every bounded subset of $G\cdot x_0$ that is metrically discrete is finite.
\end{itemize}
 \end{Lem}
\dokaz Recall that $A$ is metrically discrete if there is $s > 0$ such that $d_X(a,b) > s$
for all $a,b\in A$, $a\ne b$. Clearly, if one of Conditions 1-3 is not valid, then there is $r > 0$
such that $B(1_G,r)$ is infinite and $d_G$ is not proper. Thus, assume 1-3 hold.
 Suppose $B(1_G,r)$ is infinite for some $r > 0$ and pick $g_1$ in that set.
Suppose $\{g_n\}_{n=1}^k\subset B(1_G,2r)$ is constructed so that
$d_X(g_i\cdot x_0,x_0) < \frac{1}{i}$. Put $A=B(1_G,r)\setminus \{g_n\}_{n=1}^k$
and notice $A\cdot x_0$ is infinite (otherwise the stabilizer of $x_0$ is infinite).  Hence there are two different
elements $g,h\in A$ such that $g\cdot x_0\ne h\cdot x_0$ and $d_X(g\cdot x_0,h\cdot x_0) < \frac{1}{k+1}$. Put $g_{k+1}=g^{-1}\cdot h$.
However, $g_n\cdot x_0\to x_0$, a contradiction.
\hfill $\blacksquare$

It turns out, for nice spaces $X$, $d_G$ being a proper ls-metric is equivalent to the action being
proper.

\begin{Cor} \label{PreSvarc-Milnor}
Suppose $(X,d_X)$ is a metric space so that all infinite bounded subsets of $X$
contain an infinite Cauchy sequence.
If a group $G$ acts via isometries on $X$ and $x_0\in X$, then 
$d_G(g,h)=d_X(g\cdot x_0,h\cdot x_0)$ defines a proper left-invariant ls-metric on $G$ if and only if there is a neighborhood $U$ of $x_0$ such that the set $\{g\in G\mid g\cdot U\cap U\ne\emptyset\}$ is finite.
 \end{Cor}
\dokaz Suppose there is a neighborhood $U$ of $x_0$ such that the set $\{g\in G\mid g\cdot U\cap U\ne\emptyset\}$ is finite. Notice there is no converging sequence $g_n\cdot x_0\to x_0$
with $g_n$'s being all different.
\par If $d_G$ is proper, then choose any ball $U=B(x_0,r)$ around $x_0$.
Now, $g\cdot U\cap U\ne\emptyset$ means there is $x_g\in U$ so that $d_X(g\cdot x_g,x_0) < r$.
Therefore $d_G(g,1_G)=d_X(g\cdot x_0,x_0)\leq d_X(g\cdot x_0,g\cdot x_g)+
d_X(g\cdot x_g,x_g)+d_X(x_g,x_0)\leq r+2r+r=4r$ and there are only finitely many
such $g$'s.
\hfill $\blacksquare$

\begin{Cor} \label{BasicSvarc-Milnor}
If a group $G$ acts cocompactly and properly via isometries on a proper metric space $X$, then $g\to g\cdot x_0$ induces a coarse equivalence between $G$ and $X$
for all $x_0\in X$.
 \end{Cor}
\dokaz Define $d_G(g,h)=d_X(g\cdot x_0,h\cdot x_0)$ for all $g,h\in G$.
Clearly, $d_G$ is left-invariant. Since action is proper, $d_G$ is
a proper ls-metric. Since action is cocompact, $X$ is within bounded distance from $G\cdot x_0$.
\hfill $\blacksquare$

\begin{Cor}[\v Svarc-Milnor] \label{Svarc-Milnor}
A group
$G$ acting properly and cocompactly via isometries on a length space $X$
is finitely generated and
induces a quasi-isometry equivalence $g\to g\cdot x_0$ for any $x_0\in X$.
 \end{Cor}
\dokaz Consider the proper left-invariant metric $d_G$ induced on $G$ by the action.
The cocompactness of the action implies $G\cdot x_0$ is metrically connected.
So is $(G,d_G)$ and $G$ must be finitely generated.
Both $X$ and a Cayley graph of $G$ are proper geodesic spaces.
Therefore any coarse equivalence between them is a quasi-isometric equivalence.
\hfill $\blacksquare$

\centerline{\bf Final comments.}
Let us point out that \v Svarc-Milnor Lemma  \ref{BasicSvarc-Milnor} for non-finitely generated groups
 is useful when considering spaces of asymptotic
 dimension $0$.
A large scale analog $\MM^0$ of 0-dimensional Cantor set is
introduced in~\cite{Dran-Zar}: it is the set of all positive
integers with ternary expression containing 0's and 2's only (with
the metric from $\RR_+$): 
$$ \MM^0=\{\sum\limits_{i=-\infty}^\infty a_i3^i\mid a_i=0,2\} .$$

\begin{Prop}\cite[Theorem 3.11]{Dran-Zar}\label{M0 is universal}
The space $\MM^0$ is universal for proper metric spaces of
bounded geometry and of asymptotic dimension zero.
\end{Prop}

\begin{Prop}\label{M0 is equivalent to Z2}
The space $\MM^0$ is coarsely equivalent to
$\bigoplus\limits_{i=1}^\infty {\mathbb{Z}}_{2}$.
\end{Prop}
\dokaz Consider the subset $A=\{\sum\limits_{i=0}^\infty a_i3^i\mid a_i=0,2\}$
of $\MM^0$. Notice $ \MM^0$ is within bounded distance from $A$, so $A\to  \MM^0$ is a coarse
equivalence. Also, there is an obvious action of $\bigoplus\limits_{i=1}^\infty {\mathbb{Z}}_{2}$ on $A$ (flipping $a_i=0$ to $2$ or $a_i=2$ to $0$ if the corresponding term
in $\bigoplus\limits_{i=1}^\infty {\mathbb{Z}}_{2}$ is not zero) that is proper
and cocompact.
\hfill $\blacksquare$

Notice any infinite countable group $G$ of asymptotic dimension $0$
is locally finite (see \cite{Smith}). Thus it can be expressed as the union of
a strictly increasing sequence of its finite subgroups $G_1\subset G_2\subset\ldots$
Put $n_{1}=|G_1|$, $n_i=|G_i/G_{i-1}|$ for $i > 1$, and observe (using  
\ref{LSUniformForGroups}) that $G$ is coarsely equivalent to $\bigoplus\limits_{i=1}^\infty {\mathbb{Z}}_{n_i}$. We do not know if any two infinite countable groups of asymptotic dimension $0$ are coarsely equivalent.

\end{document}